\documentclass[11pt]{article}
\usepackage{amssymb, amsmath,  amsfonts,dsfont, mathrsfs,euscript}
\usepackage{textcomp}
\usepackage{latexsym}
\usepackage{indentfirst}
\usepackage[normalem]{ulem}

\newtheorem{Theorem}{Theorem}

\newcommand{\sgn}{\mathrm{sgn}}
\newcommand{\dd}{\,d}

\newcommand{\R}{\mathbb R}

\newcommand{\N}{\mathbb N}

\newcommand{\Realpart}{\mathrm{Re}\,}
\newcommand{\Imagpart}{\mathrm{Im}\,}

\newcommand{\ID}{\boldsymbol{I}}
\newcommand{\Q}{\boldsymbol{Q}}

\newcommand{\e}{\varepsilon}

\newcommand{\Ln}{\mathrm{Ln}\,}
\newcommand{\Arg}{\mathrm{Arg}\,}

\newcommand*\wbar[1]{%
	\hbox{%
		\vbox{%
			\hrule height 0.5pt 
			\kern0.3ex
			\hbox{%
				\kern-0.18em
				\ensuremath{#1}%
				\kern0em
			}%
		}%
	}%
}

\begin{document}
\author{A. A. Khartov$^{1,2,3,}$\footnote{Email addresses: \texttt{khartov.a@iitp.ru}, \texttt{alexeykhartov@gmail.com}} }
\title{On decomposition problem for distribution functions of class $\boldsymbol{Q}$}

\footnotetext[1]{Institute for Information Transmission Problems (Kharkevich Institute) of Russian Academy of Sciences, Bolshoy Karetny per. 19, build.1, 127051 Moscow, Russia.}
\footnotetext[2]{Laboratory for Approximation Problems of Probability, Smolensk State University, 4 Przhevalsky st., 214000 Smolensk, Russia. }
\footnotetext[3]{Saint-Petersburg National Research University of Information Technologies, Mechanics and Optics (ITMO University), 49 Kronverksky Pr., 197101 Saint-Petersburg, Russia.}

\maketitle
\begin{abstract}
	We consider a new class $\boldsymbol{Q}$ of distribution functions $F$ that have the property of rational-infinite divisibility: there exist some infinitely divisible distribution functions $F_1$ and $F_2$ such that $F_1=F*F_2$. A distribution function of the class $\boldsymbol{Q}$  is quasi-infinitely divisible in the sense that its characteristic function admits  the L\'evy--Khinchine type representation with a ``signed spectral measure''. The class $\Q$,  being a natural extension of the class $\boldsymbol{I}$ of  infinitely divisible distribution functions, is actively studied now and it finds various applications. In 2018, Lindner, Pan and Sato formulated the open question: is it true that if $F\in\boldsymbol{Q}$ and $F=F_1*F_2$ with some distribution functions $F_1$ and $F_2$, then $F_1\in\boldsymbol{Q}$ and $F_2\in\boldsymbol{Q}$? There are some  positive results under special assumptions on the type of $F$.   In this paper, we answer the question in a general setting without any additional assumptions. We also consider the same question but with the stronger assumption that $F\in\boldsymbol{I}$. 
	
\end{abstract}


\section{Introduction}

It is well known that infinitely divisible distributions play a fundamental role in probability  theory and  stochastic processes. In particular, these probability laws are the limit distributions for the sums of independent random variables in the most general shemes of summation within the uniform asymptotic negligibility of summands (see \cite{Petrov} and \cite{Zolot2}). They are also distributions of the L\'evy processes (see \cite{Sato1999}), which have further applications in the stochastic calculus   (see \cite{Appl}),  insurance mathematics (see \cite{Schoutens}) and teletraffic models (see \cite{Lifshits}).  In this article, we will consider one problem related to  the new natural extension of the class of infinitely divisible distributions.

We now turn to precise descriptions. Let $F$ be a  distribution function on the real line $\R$ with the characteristic function 
\begin{eqnarray*}
	f(t):=\int_{\R} e^{itx} \dd F(x),\quad t\in\R.
\end{eqnarray*}

We recall that $F$ (and the corresponding distribution) is called \textit{infinitely divisible} if for every positive integer $n$ there exists a  distribution function $F_{1/n}$ such that $F=(F_{1/n})^{*n}$, where ``$*$'' denotes the convolution operation, i.e. $F$ is $n$-fold convolution power of $F_{1/n}$.  It is a fundamental fact that $F$ is infinitely divisible if and only if the characteristic function $f$ is represented by \textit{the L\'evy--Khinchine formula}:
\begin{eqnarray}\label{repr_f}
	f(t)=\exp\biggl\{it \gamma+\int_{\R} \bigl(e^{itx} -1 -it \sin(x)\bigr)\tfrac{1+x^2 }{x^2} \dd G(x)\biggr\},\quad t\in\R,
\end{eqnarray}
with some \textit{shift parameter} $\gamma\in\R$, and with a  bounded non-decreasing \textit{spectral function} $G: \R \to \R$, which is assumed to be right-continuous at every point of the real line with the condition $G(-\infty)=0$. We choose $x\mapsto \sin (x)$ as the ``centering function'' in the integral in \eqref{repr_f} following to Zolotarev (see \cite{Zolot1} and \cite{Zolot2}), but  some other variants are possible (see \cite{Sato1999}, p. 38).  We note that \textit{the spectral pair} $(\gamma, G)$ is uniquely determined by $f$ (and hence by $F$). The class of all infinitely divisible distribution functions $F$ will be denoted by $\ID$. 

We call a distribution function $F$ (and the corresponding distribution) \textit{rational-infinitely divisible} if there exist some infinitely divisible distribution functions $F_1$ and $F_2$ such that $F_1=F*F_2$. Using characteristic functions, this definition is equivalent to the  admission of the formula
$f(t)=f_1(t)/f_2(t)$, $t\in\R$, where $f_1$ and $f_2$ are characteristic functions of some infinitely divisible distribution functions $F_1$ and $F_2$, respectively. Let $\Q$ denote the class of all rational-infinitely divisible distribution functions.

From the definition, it is not difficult to see that the characteristic function $f$ of $F\in\Q$ admits the \textit{L\'evy--Khinchine type representation}, i.e. \eqref{repr_f} is valid with some \textit{shift parameter} $\gamma\in\R$, and with some \textit{spectral function} $G: \R \to \R$ of bounded variation on $\R$ (in general, it is non-monototic), which is assumed to be right-continuous on $\R$ with the condition $G(-\infty)=0$. Here \textit{the spectral pair} $(\gamma, G)$ is uniquely determined by $f$  and hence by $F$ too that was already shown in \cite{GnedKolm} and \cite{Khinch}. 

Now suppose that $f$ is represented by formula \eqref{repr_f} with some $\gamma$ and  $G$ satisfying the latter conditions. Following Lindner and Sato \cite{LindPanSato}, the corresponding distribution function $F$ (and the corresponding probability law)  is called \textit{quasi-infinitely divisible}. Let us fix any real $\gamma_1$ and $\gamma_2$ such that $\gamma=\gamma_1-\gamma_2$. Due to the Hahn--Jordan decomposition, there exist some  functions $G_1$ and $G_2$, which are bounded, non-decreasing, and right-continuous on $\R$ with $G_1(-\infty)=G_2(-\infty)=0$, such that $G=G_1-G_2$. Therefore $f(t)=f_1(t)/f_2(t)$, $t\in\R$, where $f_1$ and $f_2$ are characteristic functions of infinitely divisible distribution functions with the spectral pairs $(\gamma_1, G_1)$ and $(\gamma_2, G_2)$ in canonical L\'evy--Khinchine representations for $f_1$ and $f_2$, respectively, i.e. $F$ is rational-infinitely divisible.  

Thus $F\in\Q$ if and only if $f$ admits \eqref{repr_f} with some $\gamma$ and $G$ satisfying the conditions above. So $\Q$ is an extension of the class $\ID$ of all infinitely divisible distribution functions. Namely, $F\in\ID$ if and only if $F\in\Q$ with spectral function $G$, which is non-decreasing on $\R$. Actually, the set $\Q\setminus \ID$ is non-empty and its representatives may be find in the classical monographs \cite{GnedKolm}, \cite{LinOstr}, and  \cite{Lukacs}.

The class $\Q$ is actively studied now. An extensive general analysis was done by Lindner, Pan, Sato in the paper \cite{LindPanSato}.  Criteria of belonging to this class are investigated in \cite{Alexeev_Khartov} and \cite{Khartov} for discrete probability laws, in  \cite{Berger} and \cite{BergerKutlu} for mixtures of discrete and absolutely continuous distributions, and in \cite{KhartovGenQ} for the general case. The problems on weak convergence within $\Q$ were discussed in \cite{KhartovWeak}, \cite{LindPanSato}, and  \cite{Stanek}. There are some interesting results concerning applications of the class $\Q$ (see, for instance,  \cite{ChhDemniMou}, \cite{Nakamura}, and \cite{ZhangLiuLi}).

This note is devoted to the decomposition problem of distribution functions of class $\Q$. Namely, Lindner, Pan and Sato formulated the following question (see \cite{LindPanSato}, Open Question 8.4):
\vspace{0.2cm}

\noindent\textit{Let $F$, $F_1$, $F_2$ be distribution functions on $\R$ such that $F=F_1*F_2$ and $F\in\Q$. Is it true that also $F_1$ and $F_2$  belong to the class $\Q$?}
\vspace{0.2cm}

Following the general theory (see \cite{Cuppens}, \cite{LinOstr} and \cite{Lukacs}), the equality $F=F_1*F_2$ yields a \textit{decomposition} of $F$ into the \textit{components} $F_1$ and $F_2$. The problem is rather natural and important, because if  the aswer on the question was yes, then, in particular, any component of infinitely divisible $F$ would be rational-infinitely divisible and so the class $\Q$, according to the definition of the rational-infinite divisibility, would coincide with the class of all compoponents of infinitely divisible distribution functions.  

There is a whole series of results on the problem. So the authors of \cite{LindPanSato} answered the question in the positive for the case of lattice distributions. For arbitrary discrete probability laws the similar conclusion can be easily done  based on \cite{Alexeev_Khartov} and \cite{Khartov}. Berger and Kutlu essentially generalized the previous results in the paper \cite{BergerKutlu}. Namely, they positively solved the problem for mixtures discrete and absolutely continuous distributions (with non-zero discrete part).  Also it should be noted that, in the theory of decompositions of probability laws, special representatives of $\Q$ serve as components of infinitely divisible distribution functions (see \cite{GnedKolm} p. 81--83, and \cite{LinOstr} p. 165). 

In this article, we will consider the decomposition problem for the class $\Q$ in a general setting without any additional assumptions on the type of distribution. 

We will use the following notation. We denote by $\N$  the set of positive integers. For any complex number $z$ we denote by $\Realpart\{z\}$ and $\Imagpart\{z\}$ the real and imaginary parts, respectively. The signum function is denoted by $\sgn(\cdot)$, i.e. $\sgn(x)=+1$ for $x>0$, $\sgn(x)=-1$ for $x<0$, and $\sgn(0)=0$. If $\psi$ is a complex-valued  continuous function on $\R$ satisfying $\psi(0)=1$ and $\psi(t)\ne 0$ for any $t\in\R$, then  \textit{the distinguished logarithm} $\Ln\psi$ is defined by the formula $\Ln\psi(t):=\ln |\psi(t)|+i \Arg \psi(t)$, $t\in\R$, where  $\Arg \psi(t)$ is the argument of $\psi(t)$ uniquely defined on $\R$ by the continuity with the condition $\Arg \psi(0)=0$. 

\section{Results}
Throughout this section, we construct an  example of probability law with an interesting characteristic function, which will give the answer for the formulated open question and for  one its stronger modification. We will try to describe the example rather explicitly as much as possible, but with a reasonable sufficiency. In the construction, we use the ideas from the papers  \cite{Ilinskij} and \cite{OstrFlekser}.  All conclusions will be formulated at the end of the section.

\subsection{The first auxiliary function}
Let us consider the following triangular function:
\begin{eqnarray*}
	\Delta(t):= 
	\begin{cases}
		1-|t|,& |t|\leqslant 1,\\
		0,& |t|>1.
	\end{cases}
\end{eqnarray*}
It satisfies  P\'olya's sufficient conditions (see \cite{Lukacs} p. 83) and hence it is a characteristic function of some absolutely continuous probability law. The function $\Lambda(t):=\Delta(t)e^{-|t|}$, $t\in\R$, is  characteristic too as the  product with characteristic function of  the Cauchy law. Since $\Lambda\in L_1(\R)$, it corresponds to  an absolutely continuous distribution with the following density function:
\begin{eqnarray*}
	q(x):= \dfrac{1}{2\pi}\int_{\R} e^{-itx} \Lambda(t)\dd t=\dfrac{1}{2\pi}\int_{-1}^1 e^{-itx} (1-|t|) e^{-|t|} \dd t,\quad x\in\R.
\end{eqnarray*}
In particular, we know that $q\in L_1(\R)$, $q(x)\geqslant 0$ for any $x\in\R$,  and
\begin{eqnarray}\label{eq_Lambda}
	\Lambda(t)=\int_{\R}e^{itx} q(x)\dd x,\quad t\in\R.
\end{eqnarray}

We need to find explicit expression for $q$. Observe that
\begin{eqnarray*}
	q(x)&=&\dfrac{1}{\pi}\int_{-1}^1 \cos(tx) (1-|t|)  e^{-|t|}\dd t=\dfrac{1}{\pi}\int_{-1}^0 \cos(tx)(1+t)  e^{t}\dd t\\
	&=&\dfrac{1}{2\pi}\int_{-1}^0 (1+t)\bigl(e^{(1+ix)t}+e^{(1-ix)t} \bigr)  \dd t,\quad x\in\R.
\end{eqnarray*}
We next integrate by parts with any $x\in\R$:
\begin{eqnarray*}
	q(x)&=& \dfrac{1+t}{2\pi}\biggl( \dfrac{e^{(1+ix)t}}{1+ix}+\dfrac{e^{(1-ix)t}}{1-ix}  \biggr)\bigg|_{-1}^0-\dfrac{1}{2\pi} \int_{-1}^0 \biggl( \dfrac{e^{(1+ix)t}}{1+ix}+\dfrac{e^{(1-ix)t}}{1-ix}  \biggr)\dd t\\
	&=& \dfrac{1}{2\pi}\biggl( \dfrac{1}{1+ix}+\dfrac{1}{1-ix}  \biggr)-  \dfrac{1}{2\pi} \biggl( \dfrac{1-e^{-1-ix}}{(1+ix)^2}+\dfrac{1-e^{-1+ix}}{(1-ix)^2}  \biggr)\\
	&=&\dfrac{1}{\pi(1+x^2)}-\dfrac{\ell(x)}{\pi e(1+x^2)^2},
\end{eqnarray*}
where we define
\begin{eqnarray*}
	\ell(x):=\dfrac{1}{2}\Bigl((e-e^{-ix})(1-ix)^2+(e-e^{ix})(1+ix)\Bigr),\quad x\in\R.
\end{eqnarray*}
Let us express the function $\ell$ in the real form. For any $x\in\R$ we have
\begin{eqnarray*}
	\ell(x)&=&\dfrac{1}{2}\Bigl((e-1)(1-ix)^2+(e-1)(1+ix)^2\Bigr)\\
	&&{}+\dfrac{1}{2}\Bigl((1-e^{-ix})(1-ix)^2+(1-e^{ix})\cdot(1+ix)^2\Bigr)\\
	&=&(e-1)\cdot \dfrac{1}{2}\bigl((1-ix)^2+(1+ix)^2\bigr)\\
	&&{}+\bigl(1-\cos(x)\bigr)\cdot\dfrac{1}{2}\bigl((1-ix)^2+(1+ix)^2\bigr)+\sin(x)\cdot\dfrac{i}{2}\bigl((1-ix)^2-(1+ix)^2\bigr),
\end{eqnarray*}
where
\begin{eqnarray*}
	\dfrac{1}{2}\bigl((1-ix)^2+(1+ix)^2\bigr)&=&\dfrac{1}{2}\,\bigl(2+2(ix)^2\bigr)=1-x^2,\\
	\dfrac{i}{2}\bigl((1-ix)^2-(1+ix)^2\bigr)&=&\dfrac{i}{2}\,(-4 ix)=2x.
\end{eqnarray*}
Therefore we obtain the needed formula:
\begin{eqnarray}\label{eq_ell}
	\ell(x)=(e-1) (1-x^2)+\bigl(1-\cos(x)\bigr)\cdot(1-x^2)+\sin(x)\cdot 2x,\quad x\in\R.
\end{eqnarray}

Using the latter expression for the function $\ell$, we will get lower estimates for the density $q$. So we have
\begin{eqnarray}\label{eq_q}
q(x)=\dfrac{1}{e\pi(1+x^2)}\biggl(e-\dfrac{\ell(x)}{1+x^2} \biggr),\quad x\in\R.
\end{eqnarray}
If $|x|\leqslant 1$ then
\begin{eqnarray*}
	\ell(x)\leqslant (e-1) (1-x^2)+1-\cos(x)+2|\sin(x)|\cdot |x|\leqslant (e-1) (1-x^2)+\dfrac{x^2}{2}+2 x^2.
\end{eqnarray*}
Here we applied  the well known inequalities for $\sin(x)$ and $\cos(x)$. Hence
\begin{eqnarray*}
	\ell(x)\leqslant e-1+ (3.5-e)x^2\leqslant (e-1)(1+x^2).
\end{eqnarray*}
Thus for the case $|x|\leqslant 1$ we have
\begin{eqnarray*}
	q(x)\geqslant\dfrac{1}{e\pi(1+x^2)}\biggl(e-\dfrac{(e-1)(1+x^2) }{1+x^2} \biggr)= \dfrac{1}{e\pi(1+x^2)}.
\end{eqnarray*}
If $|x|> 1$, i.e. $1-x^2<0$, then, according to \eqref{eq_ell}, we estimate
\begin{eqnarray*}
	\ell(x)\leqslant (e-1) (1-x^2)+0+2|\sin(x)|\cdot |x|\leqslant (e-1) (1-x^2)+2 |x|.	
\end{eqnarray*}
Therefore
\begin{eqnarray*}
	\ell(x)\leqslant e-1+ (1-e)x^2+2|x|= 2e-2+(1-e)(1+x^2)+2|x|.
  \end{eqnarray*}
According to \eqref{eq_q}, we have
\begin{eqnarray*}
	q(x)\geqslant\dfrac{1}{e\pi(1+x^2)}\biggl(e-\dfrac{2e-2}{1+x^2}-(1-e)-\dfrac{2|x|}{1+x^2} \biggr)=\dfrac{1}{e\pi(1+x^2)}\biggl(2e-1-\dfrac{2(e-1)}{1+x^2}-\dfrac{2|x|}{1+x^2} \biggr).
\end{eqnarray*}
Since $1+x^2>2|x|>2$ as $|x|>1$, we get
\begin{eqnarray*}
	q(x)> \dfrac{1}{e\pi(1+x^2)}\bigl(2e-1-(e-1)-1\bigr)=\dfrac{e-1}{e\pi(1+x^2)}> \dfrac{1}{e\pi(1+x^2)}.
\end{eqnarray*}

Thus we have proved that 
\begin{eqnarray}\label{ineq_q_lower_est}
	q(x)\geqslant \dfrac{1}{e\pi(1+x^2)}\quad\text{for any}\quad x\in\R.
\end{eqnarray}
By the way, it is seen from \eqref{eq_ell} that for any $x\in\R$
\begin{eqnarray*}
	|\ell(x)|\leqslant (e-1)(1+x^2)+2(1+x^2)+2|x|\leqslant (e+2)(1+x^2)
\end{eqnarray*}
(here we used the inequality $2|x|\leqslant 1+x^2$, $x\in\R$).
Therefore
\begin{eqnarray*}
	q(x)\leqslant \dfrac{1}{e\pi(1+x^2)}\biggl(e+\dfrac{|\ell(x)|}{1+x^2} \biggr)\leqslant \dfrac{2e+2}{e\pi(1+x^2)}\quad\text{for any}\quad x\in\R.
\end{eqnarray*}
This confirms that the estimate \eqref{ineq_q_lower_est} has a right order.

\subsection{The second auxiliary function}

We next fix arbitrary $n\in\N$ and we introduce a function
\begin{eqnarray*}
	\chi(t):= \exp\bigl\{-|t|+i\varphi(t)\bigr\},\quad t\in\R,
\end{eqnarray*}
where
\begin{eqnarray*}
\varphi(t):=
	\begin{cases}
		\dfrac{\bigl(t-\sgn(t)\bigr)^{2n+1}}{2n+1},& |t|>1,\\		
		0,& |t|\leqslant 1.
	\end{cases}
\end{eqnarray*}
Observe that $\varphi$ is an odd function, which is continuous on  $\R$. Moreover, in any case, $\varphi\in C^2(\R)$. Indeed,
\begin{eqnarray*}
	\varphi'(t)=
	\begin{cases}
	\bigl(t-\sgn(t)\bigr)^{2n},& |t|>1,\\		
	0,& |t|\leqslant 1,
\end{cases}
\end{eqnarray*}
and
\begin{eqnarray*}	
	\varphi''(t)=
	\begin{cases}
	2n \bigl(t-\sgn(t)\bigr)^{2n-1},& |t|>1,\\		
	0,& |t|\leqslant 1,
\end{cases}	
\end{eqnarray*}
where $2n-1\geqslant 1$.

Since $\chi\in L_1(\R)$, we can define the following function:
\begin{eqnarray}\label{def_h}
	h(x):=\dfrac{1}{2\pi} \int_{\R} e^{-itx} \chi(t) \dd t=\dfrac{1}{2\pi} \int_{\R} e^{-itx} e^{-|t|+i\varphi(t)} \dd t, \quad x\in\R.
\end{eqnarray}
We don't assert that $h$ is a density function and $\chi$ is a characteristic function of some probability law. We need only an upper estimate of $|h(x)|$ for any $x\in\R$. We first write $h(x)=J_-(x)+J_+(x)$, $x\in\R$, where
\begin{eqnarray*}
	J_-(x):=\dfrac{1}{2\pi} \int_{-\infty}^0 e^{(1-ix)t} e^{i\varphi(t)} \dd t,\qquad J_+(x):=\dfrac{1}{2\pi} \int_{0}^\infty e^{-(1+ix)t} e^{i\varphi(t)} \dd t, \quad x\in\R.
\end{eqnarray*}
Let us consider $J_-$. Integrating twice by parts, for any $x\in\R$ we get
\begin{eqnarray*}
	J_-(x)&=&\dfrac{1}{2\pi}\biggl( \dfrac{e^{(1-ix)t} }{1-ix}\cdot e^{i\varphi(t)}\biggr)\bigg|_{-\infty}^0 -\dfrac{1}{2\pi} \int_{-\infty}^0 \dfrac{e^{(1-ix)t} }{1-ix}\,i\varphi'(t) e^{i\varphi(t)} \dd t\\
	&=&\dfrac{1}{2\pi}\cdot \dfrac{1}{1-ix}-\dfrac{1}{2\pi}\biggl( \dfrac{e^{(1-ix)t} }{(1-ix)^2}\cdot i\varphi'(t) e^{i\varphi(t)}\biggr)\bigg|_{-\infty}^0\\
	&&{}+\dfrac{1}{2\pi} \int_{-\infty}^0 \dfrac{e^{(1-ix)t} }{(1-ix)^2}\,\Bigl(i\varphi''(t)+\bigl(i\varphi'(t)\bigr)^2\Bigr) e^{i\varphi(t)} \dd t\\
	&=& \dfrac{1}{2\pi(1-ix)}+\dfrac{1}{2\pi} \int_{-\infty}^1 \dfrac{e^{t-ixt} }{(1-ix)^2}\,\bigl(i\varphi''(t)-\varphi'(t)^2\bigr) e^{i\varphi(t)} \dd t.
\end{eqnarray*} 
Here we used the equalities $\varphi(t)=\varphi'(t)=\varphi''(t)=0$ for any $t\in[-1,1]$. The function $J_+$ is handled in the similar way:
\begin{eqnarray*}
	J_+(x)&=&\dfrac{1}{2\pi}\biggl( -\dfrac{e^{-(1+ix)t} }{1+ix}\cdot e^{i\varphi(t)}\biggr)\bigg|_{0}^\infty +\dfrac{1}{2\pi} \int_{0}^\infty \dfrac{e^{-(1+ix)t} }{1+ix}\,i\varphi'(t) e^{i\varphi(t)} \dd t\\
	&=&\dfrac{1}{2\pi}\cdot \dfrac{1}{1+ix}-\dfrac{1}{2\pi}\biggl( \dfrac{e^{-(1+ix)t} }{(1+ix)^2}\cdot i\varphi'(t) e^{i\varphi(t)}\biggr)\bigg|_{0}^\infty\\
	&&{}+\dfrac{1}{2\pi} \int_0^\infty \dfrac{e^{-(1+ix)t} }{(1+ix)^2}\,\Bigl(i\varphi''(t)+\bigl(i\varphi'(t)\bigr)^2\Bigr) e^{i\varphi(t)} \dd t\\
	&=& \dfrac{1}{2\pi(1+ix)}+\dfrac{1}{2\pi} \int_1^\infty \dfrac{e^{-t-ixt} }{(1+ix)^2}\,\bigl(i\varphi''(t)-\varphi'(t)^2\bigr) e^{i\varphi(t)} \dd t.
\end{eqnarray*} 
Thus for any $x\in\R$ we have
\begin{eqnarray*}
	h(x)=J_-(x)+J_+(x)=\dfrac{1}{\pi(1+x^2)}+\dfrac{1}{2\pi} \int_{\R\setminus[-1,1]} \dfrac{e^{-|t|-ixt} }{(1+ix\,\sgn (t))^2}\,\bigl(i\varphi''(t)-\varphi'(t)^2\bigr) e^{i\varphi(t)} \dd t.
\end{eqnarray*}
Therefore
\begin{eqnarray*}
	|h(x)|&\leqslant& \dfrac{1}{\pi(1+x^2)}+\dfrac{1}{2\pi} \int_{\R\setminus[-1,1]} \biggl| \dfrac{e^{-|t|-ixt} }{(1+ix\,\sgn (t))^2}\,\bigl(i\varphi''(t)-\varphi'(t)^2\bigr) e^{i\varphi(t)} \biggr|\dd t\\
	&\leqslant&\dfrac{1}{\pi(1+x^2)}+\dfrac{1}{2\pi} \int_{\R\setminus[-1,1]} \dfrac{e^{-|t|}}{1+x^2} \,\bigl(|\varphi''(t)|+|\varphi'(t)|^2\bigr)\dd t=\dfrac{e+K}{e\pi(1+x^2)},\quad x\in\R,
	\end{eqnarray*}
where
\begin{eqnarray}\label{def_K}
	K:=\dfrac{1}{2} \int_{\R\setminus[-1,1]}  \bigl(|\varphi''(t)|+|\varphi'(t)|^2\bigr)\,e^{1-|t|} \dd t.
\end{eqnarray}
Observe that for $|t|>1$
\begin{eqnarray*}
	|\varphi'(t)|= |t-\sgn(t)|^{2n}=(|t|-1)^{2n},\qquad 	|\varphi''(t)|= 2n| t-\sgn(t)|^{2n-1}=2n(|t|-1)^{2n-1}.
\end{eqnarray*}
Consequently, the integrand in \eqref{def_K} is an even function. Therefore
\begin{eqnarray*}
	K&=&\int_{1}^\infty  \bigl(2n(t-1)^{2n-1}+(t-1)^{4n}\bigr)\,e^{1-t} \dd t\\
	& =&2n\int_{0}^\infty u^{2n-1} e^{-u} \dd u +\int_{0}^\infty u^{4n}e^{-u} \dd u\\
	&=&2n\cdot (2n-1)!+(4n)!=(2n)!+(4n)!.
\end{eqnarray*}
Thus
\begin{eqnarray}\label{ineq_h}
	|h(x)|\leqslant\dfrac{e+(2n)!+(4n)!}{e\pi(1+x^2)}\quad\text{for any}\quad  x\in\R.
\end{eqnarray}
In particular, this means that $h\in L_1(\R)$ and the following formula is valid:
\begin{eqnarray}\label{eq_chi}
	\chi(t)=\int_{\R}e^{itx} h(x)\dd x,\quad t\in\R.
\end{eqnarray}

It is important to note that the function $h$ is real-valued. Indeed, returning to formula \eqref{def_h}, it is easily seen  that the integrand function
\begin{eqnarray*}
	\Lambda_x(t):=e^{-itx} e^{-|t|+i\varphi(t)}=\bigl(\cos(tx)+i\sin(tx)\bigr) \bigl( \cos\varphi(t)+i\sin\varphi(t) \bigr) e^{-|t|},\quad t\in\R,
\end{eqnarray*} 
has an odd imaginary part due to the oddness of $\varphi$. Therefore $\int_{\R}\Imagpart \bigl\{\Lambda_x(t)\bigr\}\dd t=0$ and
\begin{eqnarray*}
	h(x)=\dfrac{1}{2\pi}\int_{\R}\Lambda_x(t)\dd t=\dfrac{1}{2\pi}\int_{\R}\Realpart \bigl\{\Lambda_x(t)\bigr\}\dd t+\dfrac{i}{2\pi} \int_{\R}\Imagpart \bigl\{\Lambda_x(t)\bigr\}\dd t=\dfrac{1}{2\pi}\int_{\R}\Realpart \bigl\{\Lambda_x(t)\bigr\}\dd t\in\R.
\end{eqnarray*}

\subsection{Construction of the principal distribution}

We now fix an arbitrary $\delta>0$ such that $\delta \bigl(e+(2n)!+(4n)!\bigr)\leqslant 1$ and we define the function
\begin{eqnarray*}
	p_1(x):=\dfrac{q(x)+\delta h(x)}{1+\delta},\quad x\in\R.
\end{eqnarray*}
Since $q$ and $h$ have real values, the function $p_1$ is real-valued too. Moreover, due to \eqref{ineq_q_lower_est} and \eqref{ineq_h}, it is non-negative:
\begin{eqnarray*}
	p_1(x)\geqslant \dfrac{q(x)-\delta |h(x)|}{1+\delta}\geqslant \dfrac{1-\delta \bigl(e+(2n)!+(4n)!\bigr)}{e\pi(1+x^2)(1+\delta)}\geqslant 0,\quad x\in\R.
\end{eqnarray*}
Next,  $p_1\in L_1(\R)$, because $q$ and $h$ are from this space as we showed above. Let us define the Fourier transform of the function $p_1$:
\begin{eqnarray*}
	f_1(t):= \int_{\R} e^{itx} p_1(x)\dd x=\dfrac{1}{1+\delta}\biggl(\int_{\R} e^{itx} q(x)\dd x +\delta\int_{\R} e^{itx} h(x)\dd x\biggr),\quad t\in\R.
\end{eqnarray*}
According to \eqref{eq_Lambda} and \eqref{eq_chi}, we have
\begin{eqnarray}\label{eq_f0}
	f_1(t)=\dfrac{\Lambda(t)+\delta \chi(t)}{1+\delta},\quad t\in\R.
\end{eqnarray}
From this, by the definitions of $\Lambda$ and $\chi$, we conclude that
\begin{eqnarray*}
	\int_{\R} p_1(x)\dd x=f_1(0)=\dfrac{\Lambda(0)+\delta \chi(0)}{1+\delta}=\dfrac{1+\delta\cdot 1}{1+\delta}=1.
\end{eqnarray*}
Thus $p_1$ is a density of some probability law and $f_1$ is its characteristic function. We denote by $F_1$ the corresponding distribution function.

Let us write the explicit expression for $f_1$ using formula \eqref{eq_f0}:
\begin{eqnarray*}
	f_1(t)=\dfrac{\Delta(t)+\delta e^{i\varphi(t)}}{1+\delta}\, e^{-|t|}=
	\begin{cases}
    \tfrac{(1-|t|)+\delta}{1+\delta}\, e^{-|t|},& |t|\leqslant 1,\\
    \tfrac{\delta }{1+\delta}\,\exp\Bigl\{i\,\tfrac{(t-\sgn(t))^{2n+1}}{2n+1}\Bigr\}\, e^{-|t|},& |t|>1.
    \end{cases}
\end{eqnarray*}
From this, in particular, we find
\begin{eqnarray}\label{eq_absf_1}
		|f_1(t)|=\dfrac{\Delta(t)+\delta}{1+\delta}\, e^{-|t|},\quad t\in\R.
\end{eqnarray}
Since $\delta>0$ and $\Delta(t)\geqslant 0$ for any $t\in\R$, we conclude that $f_1(t)\ne 0$ for any $t\in\R$. Next, observe that $f_1(t)=|f_1(t)| e^{i\varphi(t)}$, $t\in\R$, and, consequently, $\Arg f_1(t)=\varphi(t)$, $t\in\R$, because $\varphi$ is continuous and $\varphi(0)=0$. So the distinguished logarithm of $f_1$ admits the following representation
\begin{eqnarray*}
	\Ln f_1(t)&=&\ln|f_1(t)|+i\Arg f_1(t)\\
	&=&\ln\biggl( \dfrac{\Delta(t)+\delta}{1+\delta} \biggr)-|t|+i \varphi(t)
	=\begin{cases}
		\ln\bigl(\tfrac{(1-|t|)+\delta}{1+\delta}\bigr) -|t|,& |t|\leqslant 1,\\
		\ln\bigl(\tfrac{\delta }{1+\delta}\bigr)-|t|+i\,\tfrac{(t-\sgn(t))^{2n+1}}{2n+1},& |t|>1.
	\end{cases}
\end{eqnarray*}

We now show that $F_1\notin \Q$. Let us write explicitly the values of  $\Ln f_1(t)$ for $t>1$:
\begin{eqnarray*}
	\Ln f_1(t)=\ln\biggl(\dfrac{\delta }{1+\delta}\biggr)-t+i\,\dfrac{(t-1)^{2n+1}}{2n+1}.
\end{eqnarray*}
This, in particular, implies that
 \begin{eqnarray*}
 	\lim_{t\to\infty}\dfrac{|\Ln f_1(t)|}{t^{2n+1}}= \dfrac{1}{2n+1},
 \end{eqnarray*}
where $2n+1\geqslant 3$. Therefore $|\Ln f_1(t)|/t^2\to \infty$ as $t\to\infty$, which is impossible for characteristic functions of (rational-)infinitely divisible laws due to the following fact: if $f$ is a characteristic function of some $F\in\Q$ (or $\ID$), then $\Ln f(t)/t^2\to c$ as $t\to\infty$ with some constant $c\leqslant 0$ (see \cite{LindPanSato}, Lemma 2.7, or \cite{Sato1999}, Lemma 43.11). Thus we proved that $F_1\notin \Q$.

Let us consider the function $|f_1|(t):=|f_1(t)|$, $t\in\R$, which is described by formula \eqref{eq_absf_1}. It is an interesting fact that $|f_1|$ is a characteristic function of some distribution function $F_{1,a}\in\Q$. We will prove it and find the corresponding spectral pair. 

Recall that $t\mapsto e^{-|t|}$, $t\in\R$, is the characteristic function of distribution function $F_C$ of the Cauchy law, which is infinitely divisible, i.e. $F_C\in\ID$. The corresponding L\'evy-Khinchine representation is well known (see \cite{Zolot1} pp. 4--8):
\begin{eqnarray*}
	e^{-|t|}=  \exp\biggl\{\int_{\R\setminus\{0\}}\bigl(e^{itx}-1-it\sin(x)\bigr)    \dfrac{\dd x}{\pi x^2}    \biggr\}=\exp\biggl\{\int_{\R}\bigl(e^{itx}-1-it\sin(x)\bigr)\tfrac{1+x^2}{x^2}\dd G_{C}(x)\biggr\},\quad t\in\R,
\end{eqnarray*}
where the shift parameter  equals zero and the spectral function $G_C$ is defined by the formula:
\begin{eqnarray*}
	G_{C}(x):=\int_{-\infty}^{x}\dfrac{\dd u}{\pi(1+u^2)},\quad x\in\R.
\end{eqnarray*}

Next, we mentioned above that $t\mapsto\Delta(t)$, $t\in\R$,  is a characteristic function of some absolutely continuous distribution. Therefore  $t\mapsto \tfrac{\Delta(t)+\delta}{1+\delta}$, $t\in\R$ is a characteristic function of some distribution, which is a  mixture of the previous one and  the degenerate law concentrated at the point $x=0$. Let $F_m$ denote its distribution function. Let us show that $F_m\in\Q$ and find the corresponding spectral pair. We write
\begin{eqnarray}\label{eq_MixDelta_exp}
	\dfrac{\Delta(t)+\delta}{1+\delta}=\exp\biggl\{ \ln\biggl(\dfrac{\delta+\Delta(t)}{1+\delta} \biggr)   \biggr\}=\exp\biggl\{ \ln\biggl(1+\dfrac{\Delta(t)}{\delta} \biggr) -\ln\biggl(1+\dfrac{1}{\delta} \biggr)  \biggr\},\quad t\in\R,
\end{eqnarray}
where the function $t\mapsto  \ln \bigl(1+\tfrac{1}{\delta}\Delta(t)\bigr)$, $t\in\R$, is continuous on $\R$ and it has the support $[-1,1]$. So it belongs to $L_1(\R)$ and we can define the function
\begin{eqnarray*}
	g_m(x):=\dfrac{1}{2\pi}\int_{\R}e^{-itx}  \ln\biggl(1+\dfrac{\Delta(t)}{\delta}\biggr)\dd t,\quad x\in\R.
\end{eqnarray*}
By the definition of $\Delta(t)$, we have
\begin{eqnarray}\label{eq_gm1}
	g_m(x)=\dfrac{1}{2\pi}\int_{-1}^1 e^{-itx} \ln\biggl( 1+\dfrac{1-|t|}{\delta} \biggr)\dd t=\dfrac{1}{\pi}\int_{0}^1 \cos(tx) \ln\biggl( 1+\dfrac{1-t}{\delta} \biggr)\dd t,\quad x\in\R.
\end{eqnarray}
Observe that for any $x\in\R$
\begin{eqnarray}\label{ineq_abs_gm}
	|g_m(x)|\leqslant \dfrac{1}{\pi}\int_{0}^1 \ln\biggl( 1+\dfrac{1-t}{\delta} \biggr)\dd t\leqslant \dfrac{1}{\pi}\int_{0}^1 \dfrac{1-t}{\delta} \dd t=\dfrac{1}{2\pi\delta}.
\end{eqnarray}
For $x\ne 0$ we apply the  integration by parts in \eqref{eq_gm1}:
\begin{eqnarray*}
	g_m(x)=\biggl(\dfrac{\sin(tx)}{\pi x}\, \ln\biggl( 1+\dfrac{1-t}{\delta} \biggr)\biggr)\bigg|_0^1-\dfrac{1}{\pi }\int_{0}^1 \dfrac{\sin(tx)}{x}\cdot \dfrac{-1}{\delta+1-t}\dd t=\dfrac{1}{\pi x}\int_{0}^1  \dfrac{\sin(tx)}{\delta+1-t}\dd t.
\end{eqnarray*}
We again integrate by parts :
\begin{eqnarray*}
	g_m(x)&=&-\biggl(\dfrac{1}{\pi x}\cdot\dfrac{\cos(tx)}{x}\cdot \dfrac{1}{\delta+1-t}\biggr)\bigg|_0^1+\dfrac{1}{\pi x}\int_{0}^1 \dfrac{\cos(tx)}{x}\, \dfrac{1}{(\delta+1-t)^2}\dd t\\
	&=&\dfrac{1}{\pi x^2(1+\delta)}-\dfrac{\cos(x)}{\pi x^2\delta}+\dfrac{1}{\pi x^2}\int_{0}^1  \dfrac{\cos(tx)}{(\delta+1-t)^2}\dd t.
\end{eqnarray*}
Thus we obtain that
\begin{eqnarray}\label{eq_gm2}
	g_m(x)=\dfrac{1}{\pi\delta x^2}\,\biggl(\dfrac{\delta}{1+\delta}-\cos(x)+\delta\e(x)\biggr),\quad x\ne 0,
\end{eqnarray}
where
\begin{eqnarray*}
	\e(x):= \int_{0}^1  \dfrac{\cos(tx)}{(\delta+1-t)^2}\dd t,\quad x\ne 0.
\end{eqnarray*}
Since the function $t\mapsto \tfrac{1}{(\delta+1-t)^2}$ is continuous and bounded on $[0,1]$, the function $\e$ is bounded on $\R$. From the classical Fourier analysis we also know that $\e(x)\to 0$ as $x\to \pm\infty$. So, due to \eqref{ineq_abs_gm} and \eqref{eq_gm2}, we conclude that $g_m\in L_1(\R)$.  Hence the following inversion formula holds
\begin{eqnarray*}
	\ln\biggl(1+\dfrac{\Delta(t)}{\delta}\biggr)=\int_{\R}e^{itx}  g_m(x)\dd x,\quad t\in\R.
\end{eqnarray*}
We note that
\begin{eqnarray*}
	\ln\biggl(1+\dfrac{1}{\delta}\biggr)=\ln\biggl(1+\dfrac{\Delta(0)}{\delta}\biggr)=\int_{\R}  g_m(x)\dd x.
\end{eqnarray*}
Returning to \eqref{eq_MixDelta_exp}, we write
\begin{eqnarray*}
	\dfrac{\Delta(t)+\delta}{1+\delta}=\exp\biggl\{\int_{\R}(e^{itx}-1)  g_m(x)\dd x   \biggr\}=\exp\biggl\{\int_{\R}\bigl(e^{itx}-1-it\sin(x)\bigr)  g_m(x)\dd x   \biggr\},\quad t\in\R,
\end{eqnarray*}
where $it \sin(x)$ were formally added due to the evenness of the function $g_m$ (see formula \eqref{eq_gm1}). Thus we come to the L\'evy--Khinchine type representation
\begin{eqnarray*}
	\dfrac{\Delta(t)+\delta}{1+\delta}=\exp\biggl\{\int_{\R}\bigl(e^{itx}-1-it\sin(x)\bigr)\tfrac{1+x^2}{x^2} \dd G_m(x)  \biggr\},\quad t\in\R,
\end{eqnarray*}
with zero shift parameter and the spectral function
\begin{eqnarray*}
	G_m(x):=\int_{-\infty}^x \dfrac{u^2g_m(u)}{1+u^2}\,\dd u,\quad x\in\R. 
\end{eqnarray*}
Therefore $F_m$ is quasi-infinitely divisible and, consequently, $F_m\in\Q$ (see the  introduction). By the way, the very fact that $F_m\in\Q$ follows by Berger's criterion from \cite{Berger}. We didn't apply it because we need the explicit expression for the spectral function.

We now return to $|f_1|$. According to \eqref{eq_absf_1}, we see that $|f_1|$ is the  characteristic function corresponding to the distribution function $F_{1,a}:=F_m*F_C$ and $|f_1|$ admits the L\'evy-Khinchine type representation:
\begin{eqnarray}\label{repr_f1}
	|f_1(t)|=\dfrac{\Delta(t)+\delta}{1+\delta}\, e^{-|t|}=\exp\biggl\{\int_{\R}\bigl(e^{itx}-1-it\sin(x)\bigr)\tfrac{1+x^2}{x^2} \dd G_{1,a}(x)  \biggr\},\quad t\in\R,
\end{eqnarray}
with zero shift parameter and the spectral function  $G_{1,a}(x):=G_m(x)+G_C(x)$, $x\in\R$, i.e.
\begin{eqnarray*}
	G_{1,a}(x)=\int_{-\infty}^x \dfrac{u^2g_m(u)}{1+u^2}\,\dd u +\int_{-\infty}^x \dfrac{1}{\pi(1+u^2)}\,\dd u=\int_{-\infty}^x \dfrac{\pi u^2g_m(u)+1}{\pi(1+u^2)}\,\dd u,\quad x\in\R.
\end{eqnarray*}
Therefore $F_{1,a}\in \Q$. Actually, we can prove more that $F_{1,a}\in\Q\setminus \ID$.  Indeed, according to \eqref{eq_gm2}, for any $u\ne 0$
\begin{eqnarray*}
	\pi u^2g_m(u)+1=\dfrac{1}{\delta }\,\biggl(\dfrac{\delta}{1+\delta}-\cos(u)+\delta\e(u)\biggr)+1=\dfrac{1}{\delta }\,\biggl(\delta+\dfrac{\delta}{1+\delta}-\cos(u)+\delta\e(u)\biggr).
\end{eqnarray*}
Hence  the inequality $\pi u^2g_m(u)+1<\tfrac{1}{\delta }\,\bigl(2\delta-\cos(u)\bigr)$ holds for all $u$ with sufficiently large absolutely values and, in particular, for these $u$ we have $\pi u^2g_m(u)+1<0$ on the intervals, where $\cos(u)>2\delta$. Therefore $G_{1,a}$  cannot be monotonic. This means that $F_{1,a}\in\Q\setminus \ID$. By the way, as corollary, $F_{m}\in\Q\setminus \ID$ too.

\subsection{The statements}

We now introduce a distribution function $F:=F_{1,a}*F_{1,a}$. Let $f$ denote its characteristic function. On the one hand, we have that $f(t)=|f_1(t)|^2$, $t\in\R$. On account of \eqref{repr_f1}, the function $f$ admits the L\'evy-Khinchine type representation \eqref{repr_f} with $\gamma:=0$ and $G(x):=2G_{1,a}(x)$, $x\in\R$, and, consequently, $F\in\Q$ and even more $F\in\Q\setminus \ID$. On the other hand, $f(t)=|f_1(t)|^2=f_1(t)\overline{f_1(t)}$, $t\in\R$, which means that $F=F_1*\,\wbar{F}_1$, where $\wbar{F}_1$ denotes the conjugate distribution function of $F_1$ (see \cite{Lukacs}, p. 29), i.e. $\wbar{F}_1(x):= 1-F_1(-x-0)$, $x\in\R$. It is clear that $\wbar{F}_1\notin \Q$ together with $F_1$ (see description of $\Q$ in the  introduction).  

Thus we come to the following assertion.
\begin{Theorem}
	There exist distribution functions $F$, $F_1$, and $F_2$ such that $F=F_1*F_2$ and $F\in\Q\setminus \ID$, but $F_1\notin \Q$ and $F_2\notin \Q$.  
\end{Theorem} 
So we answer the open question from the introduction  in the negative. 

Will the answer change in the positive if we modify the formulation of the question replacing the assumption that $F\in\Q$ by the stronger one that $F\in\ID$? In order to answer this question, we continue to work with our example within the same notations. Recall that $F\in\Q$ with the spectral pair $(0,G)$ in the representation of its characteristic function $f$, where  
\begin{eqnarray*}
	G(x)=2\int_{-\infty}^x \dfrac{\pi u^2g_m(u)+1}{\pi(1+u^2)}\,\dd u,\quad x\in\R.
\end{eqnarray*} 
Let us define non-negative functions
\begin{eqnarray*}
	g_+(u):=\max\bigl\{0,\pi u^2g_m(u)+1\bigr\}\quad\text{and}\quad g_-(u):=-\min\bigl\{0,\pi u^2g_m(u)+1\bigr\},\quad u\in\R.
\end{eqnarray*}
So we have $\pi u^2g_m(u)+1=g_+(u)-g_-(u)$, $u\in\R$.  Then we write $G(x)=G_+(x)-G_-(x)$, $x\in\R$, where
\begin{eqnarray*}
	G_+(x)=\int_{-\infty}^x \dfrac{2g_+(u)}{\pi(1+u^2)}\,\dd u\quad\text{and}\quad G_-(x)=\int_{-\infty}^x \dfrac{2g_-(u)}{\pi(1+u^2)}\,\dd u,\quad x\in\R.
\end{eqnarray*}
The functions $G_+$ and $G_-$ are non-decreasing on $\R$ due to the non-negativity of $g_+$ and $g_-$. Let $F_+$ and $F_-$ be infinitely divisible distribution functions with the spectral pairs $(0,G_+)$ and $(0,G_-)$, respectively. Their characteristic functions, say $f_+$ and  $f_-$, respectively, satisfy the equality $f(t)=f_+(t)/f_-(t)$, $t\in\R$.
Writing it in the  form $f_+(t)=f(t)f_-(t)$, $t\in\R$, we conclude that $F_+=F*F_-$. Then $F_+=(F_1*\,\wbar{F}_1)*F_-=F_1*(\,\wbar{F}_1*F_-)$, where $F_+\in\ID$ and $F_1\notin \Q$. Actually,  $(\,\wbar{F}_1*F_-)\notin \Q$, since otherwise we can find infinitely divisible distribution function $F_3$ and $F_4$ such that $F_3= (\,\wbar{F}_1*F_-)*F_4$ and, consequently,
\begin{eqnarray*}
	F_+*F_4=F_1*(\,\wbar{F}_1*F_-)*F_4=F_1*F_3,
\end{eqnarray*}
where $(F_+*F_4)$ and $F_3$ are infinitely divisible, i.e. $F_1\in\Q$, a contradiction. 
\begin{Theorem}
	There exist distribution functions $F$, $F_1$, and $F_2$ such that $F=F_1*F_2$ and $F\in\ID$, but $F_1\notin \Q$ and $F_2\notin \Q$.  
\end{Theorem} 
So we answer the modified open question in the negative too.

The reason of the negative results, as we saw above, is that a component of a (rational-)infinitely divisible distribution function can have the characteristic function, whose argument has a very fast growth at the infinity. In our example, the argument of the component characteristic function $f_1$ is a power function of arbitrarily high degree.

Thus Theorem 1 and 2 show that in order to get the positive solution of the decomposition problem for the classes $\Q$ and $\ID$ we need to make some additional assumptions for the distribution function $F$.

\end{document}